*Original Article*

# Factors of Composite $4n^2 + 1$ using Fermat's Factorization Method

Paul Ryan A. Longhas[1], Alsafat M. Abdul[2], Aurea Z. Rosal[3]

*[1,2]Instructor, Department of Mathematics and Statistics,*
*Polytechnic University of the Philippines, Manila, Philippines*

*[3]Associate Professor, Department of Mathematics and Statistics,*
*Polytechnic University of the Philippines, Manila, Philippines*

[1]pralonghas@pup.edu.ph, [2]amabdul@pup.edu.ph, [3]azrosal@pup.edu.ph

*Abstract – In this article, we factor the composite $4n^2 + 1$ using Fermat's factorization method. Consequently, we characterized all proper factors of composite $4n^2 + 1$ in terms of its form. Furthermore, the composite Fermat's number is considered in this study.*

## I. INTRODUCTION

Fermat's factorization is a method in factoring the odd composite natural number $N$ by expressing it as difference of two squares. If $N = ab$ where $N$ is odd positive composite number, and $a$ and $b$ are proper factor of $N$, then $N$ can be written as
$$N = c^2 - d^2 = (c+d)(c-d)$$
where $c = \frac{a+b}{2}$ and $d = \frac{b-a}{2}$[1]. For example, if we want to factor 9797 using Fermat's factorization method, then the goal is to expressed 9797 as a difference of two square, that is, $9797 = c^2 - d^2$ for some $c, d \in \mathbb{Z}$. In solving $c$ and $d$, note that $d^2 = c^2 - 9797$ implies that $c \geq \lceil \sqrt{9797} \rceil$, that is, $c \geq 99$. If $c = 99$, then $d^2 = 99^2 - 9797 = 4$ is a perfect square, and thus, $9797 = (99-2)(99+2)$. In general, the Fermat's method might be slower than trial and error method to apply. In fact, Fermat's factorization works best to $N$ when there is a factor $a$ of $N$ such that $a$ is near to $\sqrt{N}$ [1]. Thus, some improvement is necessary to make the Fermat's method effective [1,7,8,9]. In 1999, R. Lehman devised a systematic method to improve the Fermat's method by multiplier improvement so that the Fermat's method plus trial division can be factor $N$ in $O\left(N^{\frac{1}{3}}\right)$ time [1].

In this article, we study the Fermat's factorization of composite $4n^2 + 1$ and its proper factors. More precisely speaking, we proved that:

1. If $n$ is even, then every proper factors of composite $N = 4n^2 + 1$ is can be expressed as
$$8u + 1 \pm \sqrt{(8u+1)^2 - N}$$
where $u \in \mathbb{N}$ and
1.1. $(8u+1) - N$ is a perfect square.
1.2. $u \in \left[\frac{-1+\sqrt{N}}{8}, \frac{N-5}{40}\right)$.
1.3. $u \not\equiv 0 \pmod{p}$ for all prime $p \equiv 3 \pmod{4}$.
1.4. For all odd prime $p$ does not divide $N$, $u \equiv \frac{x_0^{-1}N + x_0 - 2}{16} \pmod{p}$, for some $x_0 \in \mathbb{Z}_p \setminus \{0\}$.

2. If $n$ is odd, then every proper factors of composite $N = 4n^2 + 1$ is can be expressed as
$$8u + 3 \pm \sqrt{(8u+3)^2 - N}$$
where $u \in \mathbb{N}$ and
2.1. $(8u+3)^2 - N$ is a perfect square.

2.2. $u \in \left[\frac{-3+\sqrt{N}}{8}, \frac{N-15}{40}\right)$.

2.3. $4u + 1 \not\equiv 0 \pmod{p}$ for all prime $p \equiv 3 \pmod{4}$.

2.4. For all odd prime $p$ does not divide $N$, $u \equiv \frac{x_0^{-1}N + x_0 - 6}{16} \pmod{p}$, for some $x_0 \in \mathbb{Z}_p \setminus \{0\}$.

The main results characterized all proper factors of $4n^2 + 1$ in terms of its form.

## II. MAIN RESULTS

First, we study $4n^2 + 1$ when $n$ is even. The following lemma is vital in the proof of main results of this study. The proof of Lemma 2.1 follows from the fact that every factor of $16m^2 + 1$ is can be expressed as $4a + 1$ where $a$ is a positive integer [5].

**Lemma 2.1.** Let $m \in \mathbb{N}$. If $16m^2 + 1$ is composite, then there is a natural number $b \leq \frac{-1+\sqrt{16m^2+1}}{4}$ where

$$m^2 + b^2 \equiv 0 \pmod{(4b+1)}. \tag{1}$$

Furthermore, $4b + 1$ is a proper factor of $16m^2 + 1$.

**Proof:** Assume $16m^2 + 1$ is composite. Then, by [5] there exists a natural number $a$ and $b$ such that

$$16m^2 + 1 = (4a+1)(4b+1) \tag{2}$$

where $4b + 1 \leq \sqrt{16m^2 + 1}$ is a proper factor of $16m^2 + 1$. Manipulate equation in (2), then we have

$$4m^2 = 4ab + a + b. \tag{3}$$

So, there exists $u \in \mathbb{N}$ such that $4u = a + b$. Replacing $a = 4u - b$ and $a + b = 4u$ in equation (3), then we obtain

$$4m^2 = 4(4u - b)b + 4u \tag{4}$$

which gives

$$m^2 + b^2 = u(4b + 1) \tag{5}$$

where $b \leq \frac{-1+\sqrt{16m^2+1}}{4}$, as desired. □

The following proposition follows from lemma 2.1.

**Proposition 2.2.** Let $m \in \mathbb{N}$. If $N = 16m^2 + 1$ is composite, then there exists a natural number $u$ such that

$$N = \left(8u + 1 + \sqrt{(8u+1)^2 - N}\right)\left(8u + 1 - \sqrt{(8u+1)^2 - N}\right) \tag{6}$$

where

1. $(8u + 1) - N$ is a perfect square.

2. $u \in \left[\frac{-1+\sqrt{N}}{8}, \frac{N-5}{40}\right)$.

3. $u \not\equiv 0 \pmod{p}$ for all prime $p \equiv 3 \pmod{4}$.

4. For all odd prime $p$ does not divide $N$, $u \equiv \frac{x_0^{-1}N + x_0 - 2}{16} \pmod{p}$, for some $x_0 \in \mathbb{Z}_p \setminus \{0\}$.

Furthermore, the following holds:

5. If $n$ is even, then $u \equiv 2 \pmod{4}$.

6. If $n$ is not divisible by 3, then $u \equiv 1 \pmod{3}$.

Conversely, if there exists $u \in \left[\frac{-1+\sqrt{N}}{8}, \frac{N-5}{40}\right)$ where $(8u + 1) - N$ is a perfect square, then $N$ is composite.

**Proof:** Assume $N = 16m^2 + 1$ is composite, then by lemma 2.1 there exists a natural number $b$ and natural number $u$ where

$$u = \frac{m^2 + b^2}{4b + 1} \tag{7}$$

and $4b + 1$ is a proper factor of $N$. Consider the quadratic polynomial defined by
$$f(x) = x^2 - 2(8u + 1)x + N. \tag{8}$$
Note that $f(4b + 1) = 0$ and the product of the roots of $f(x)$ is $N$, and thus, every roots of $f(x)$ is a proper factor of $N$. Computing the zeroes $r_1$ and $r_2$ of $f(x)$ yields
$$r_1 = 8u + 1 - \sqrt{(8u + 1)^2 - N} \tag{9}$$
$$r_2 = 8u + 1 + \sqrt{(8u + 1)^2 - N}. \tag{10}$$
Thus,
$$N = \left(8u + 1 + \sqrt{(8u + 1)^2 - N}\right)\left(8u + 1 - \sqrt{(8u + 1)^2 - N}\right). \tag{11}$$

1. Since $r_1 = 8u + 1 - \sqrt{(8u + 1)^2 - N} \in \mathbb{N}$, then $(8u + 1) - N$ is a perfect square.

2. Note that we have
$$(8u + 1)^2 - N \geq 0. \tag{12}$$
Solve the inequality in (12) with the assumption that $u > 0$ yields
$$u \geq \frac{-1 + \sqrt{N}}{8}. \tag{13}$$
Furthermore, by Rolle's theorem, there exists $\theta \in (r_1, r_2) \subset \left(r_1, \frac{N}{5}\right]$ where $f'(\theta) = 0$. Thus,
$$2\theta - 2(8u + 1) = 0. \tag{14}$$
Hence, $u = \frac{\theta - 1}{8}$. Since $\theta \in (r_1, r_2) \subset \left(r_1, \frac{N}{5}\right]$, then $u \leq \frac{N-5}{40}$. Therefore, we will have
$$u \in \left[\frac{-1 + \sqrt{N}}{8}, \frac{N - 5}{40}\right). \tag{15}$$

3. We claim that $u \not\equiv 0 \pmod{p}$ for all prime $p \equiv 3 \pmod{4}$. Indeed, assume there is a prime $p \equiv 3 \pmod{4}$ where $u \equiv 0 \pmod{p}$. Then, we have
$$r_1^2 - 2(8u + 1)r_1 + N = 0$$
$$r_1^2 - 2(8u + 1)r_1 + N \equiv 0 \pmod{p}$$
$$r_1^2 - 2r_1 + 16m^2 \equiv 0 \pmod{p}$$
$$(r_1 - 1)^2 + (4m)^2 \equiv 0 \pmod{p}.$$
Note that $(r_1 - 1)^2 + (4m)^2 \equiv 0 \pmod{p}$ is impossible for $p \equiv 3 \pmod{4}$, a contradiction.

4. Let $p$ be an odd prime where $N \not\equiv 0 \pmod{p}$. Then, we claim that $u \equiv \frac{x_0^{-1}N + x_0 - 2}{16} \pmod{p}$, for some $x_0 \in \mathbb{Z}_p \setminus \{0\}$. Indeed, since $f(4b + 1) = 0$, then $f(x)$ is reducible over $\mathbb{Z}$, and thus, $f(x)$ is reducible over $\mathbb{Z}_p$. Therefore, there is $x_0 \in \mathbb{Z}_p$ such that
$$x_0^2 - 2(8u + 1)x_0 + N \equiv 0 \pmod{p}. \tag{16}$$
Since $N \not\equiv 0 \pmod{p}$, then $x_0 \neq 0$. Thus, from equation in (16) we have
$$u \equiv \frac{x_0^{-1}N + x_0 - 2}{16} \pmod{p} \tag{17}$$
for some $x_0 \in \mathbb{Z}_p \setminus \{0\}$.

5. Suppose $n$ is even. We claim that $u \not\equiv 2 \pmod{4}$. Indeed, if $u \equiv 2 \pmod{4}$, then there exists $k \in \mathbb{Z}$ where $u = 4k + 2$. Thus, we have
$$f(x) = x^2 - 2(8(4k + 2) + 1)x + N. \tag{18}$$
Replacing $x = 4y + 1$ in equation $x^2 - 2(8(4k + 2) + 1)x + N = 0$, then we have
$$y^2 - 4(4k + 2)y - (4k + 2) + n^2 = 0. \tag{19}$$
Since $f(x)$ is reducible over $\mathbb{Z}$, then the polynomial $g(x) = x^2 - 4(4k + 2)x - (4k + 2) + n^2$ is also reducible over $\mathbb{Z}$ in which one of the zeroes is $b$. Notice that $2 \nmid 1, 2|4(4k + 2), 2|(-(4k + 2) + n^2)$ but $2^2 \nmid (-(4k + 2) + n^2)$, so by Eisenstein Criterion theorem [6], $g(x)$ is irreducible over $\mathbb{Q}$ which is a contradiction. Therefore, $u \not\equiv 2 \pmod{4}$.

6. The proof is follow from statement 3 by setting $p = 3$.

Conversely, if there exists $u \in \left[\frac{-1+\sqrt{N}}{8}, \frac{N-5}{40}\right)$ where $(8u + 1) - N$ is a perfect square, then
$$N = \left(8u + 1 + \sqrt{(8u + 1)^2 - N}\right)\left(8u + 1 - \sqrt{(8u + 1)^2 - N}\right). \tag{20}$$
If $N$ is prime, then $8u + 1 - \sqrt{(8u + 1)^2 - N} = 1$, and hence, $u = m^2$ which is a contradiction since $u \in \left[\frac{-1+\sqrt{N}}{8}, \frac{N-5}{40}\right)$. Therefore, $N$ is composite. □

The equation in (6) is called the Fermat's factorization of $16m^2 + 1$. Note that statements 1-6 of proposition 2.2 states the property of all factors of $16m^2 + 1$ if we factor $16m^2 + 1$ using Fermat's factorization method. This also gives a sieve method to determine the proper factor of composite $16m^2 + 1$. The next result characterized all proper factor of composite $16m^2 + 1$ in terms of its structure.

**Proposition 2.3.** Every proper factors of composite $N = 16m^2 + 1$ is can be expressed as
$$8u + 1 \pm \sqrt{(8u + 1)^2 - N} \tag{21}$$
where $u \in \mathbb{N}$ and
1. $(8u + 1) - N$ is a perfect square.
2. $u \in \left[\frac{-1+\sqrt{N}}{8}, \frac{N-5}{40}\right)$.
3. $u \not\equiv 0 \pmod{p}$ for all prime $p \equiv 3 \pmod{4}$.
4. For all odd prime $p$ does not divide $N$, $u \equiv \frac{x_0^{-1}N + x_0 - 2}{16} \pmod{p}$, for some $x_0 \in \mathbb{Z}_p \setminus \{0\}$.

Furthermore, the following holds:
5. If $n$ is even, then $u \equiv 2 \pmod{4}$.
6. If $n$ is not divisible by 3, then $u \equiv 1 \pmod{3}$.

**Proof:** The results follow directly from Proposition 2.2. Furthermore, (21) follows from (9) and (10). □

Now, we study the factorization of $4n^2 + 1$ where $n$ is odd number using Fermat's factorization method.

**Proposition 2.4.** Let $m \in \mathbb{N}$. If $N = 4(2m + 1)^2 + 1$ is composite, then there is $u \in \mathbb{N}$ such that
$$N = \left(8u + 3 + \sqrt{(8u + 3)^2 - N}\right)\left(8u + 3 - \sqrt{(8u + 3)^2 - N}\right). \tag{22}$$
where
1. $(8u + 3)^2 - N$ is a perfect square.
2. $u \in \left[\frac{-3+\sqrt{N}}{8}, \frac{N-15}{40}\right)$.
3. $4u + 1 \not\equiv 0 \pmod{p}$ for all prime $p \equiv 3 \pmod{4}$.
4. For all odd prime $p$ does not divide $N$, $u \equiv \frac{x_0^{-1}N + x_0 - 6}{16} \pmod{p}$, for some $x_0 \in \mathbb{Z}_p \setminus \{0\}$.

Conversely, if there exists $u \in \left[\frac{-3+\sqrt{N}}{8}, \frac{N-15}{40}\right)$ where $(8u + 3) - N$ is a perfect square, then $N$ is composite.

**Proof:** Let $a, b \in \mathbb{N}$ such that $4(2m + 1)^2 + 1 = (4a + 1)(4b + 1)$. Then, $a + b = 4(m^2 + m - ab) + 1 \in \mathbb{N}$. Take $u = m^2 + m - ab \in \mathbb{N}$. Then,
$$a + b = 4u + 1. \tag{23}$$

Applying Fermat's factorization method in $4(2m+1)^2 + 1$, then we have
$$N = \left(2a + 2b + 1 + \sqrt{(2a+2b+1)^2 - N}\right)\left(2a + 2b + 1 - \sqrt{(2a+2b+1)^2 - N}\right). \tag{24}$$
Thus, we have
$$N = \left(8u + 3 + \sqrt{(8u+3)^2 - N}\right)\left(8u + 3 - \sqrt{(8u+3)^2 - N}\right). \tag{25}$$

1. Since $8u + 3 + \sqrt{(8u+3)^2 - N} \in \mathbb{N}$, then $(8u+3)^2 - N$ is a perfect square.

2. Note that we have
$$(8u+3)^2 - N \geq 0. \tag{26}$$
Solving the inequality in (26) with the assumption that $u > 0$, then we have
$$u \geq \frac{-3 + \sqrt{N}}{8}. \tag{27}$$
Consider the quadratic function
$$f(x) = x^2 - 2(8u+3)x + N. \tag{28}$$
Since $f(4b+1) = 0$ and the product of the roots of $f(x)$ is $N$, then all roots of $f$ is a proper factor of $N$. Let $r_1$ and $r_2$ be the roots of $f(x)$. Applying Rolle's theorem, then there is $\theta \in (r_1, r_2) \subset \left(r_1, \frac{N}{5}\right]$ where $f'(\theta) = 0$. Thus, we have
$$2\theta - 2(8u+3) = 0. \tag{29}$$
So, $\theta = 8u + 3 > \frac{N}{5}$, and thus, $u > \frac{N-15}{40}$. Therefore,
$$u \in \left[\frac{-3 + \sqrt{N}}{8}, \frac{N-15}{40}\right). \tag{30}$$

3. We claim that $4u + 1 \not\equiv 0 \pmod{p}$ for all prime $p \equiv 3 \pmod{4}$. Indeed, assume there is a prime $p \equiv 3 \pmod{4}$ where $4u + 1 \equiv 0 \pmod{p}$. Then, we have
$$r_1^2 - 2(8u+3)r_1 + N = 0$$
$$r_1^2 - 2(2(4u+1)+1)r_1 + N \equiv 0 \pmod{p}$$
$$r_1^2 - 2r_1 + 16m^2 \equiv 0 \pmod{p}$$
$$(r_1 - 1)^2 + (4m)^2 \equiv 0 \pmod{p}.$$
Note that $(r_1 - 1)^2 + (4m)^2 \equiv 0 \pmod{p}$ is impossible for $p \equiv 3 \pmod{4}$, a contradiction.

4. Since $f$ is reducible over $\mathbb{Z}$, then $f$ is reducible over $\mathbb{Z}_p$, for all odd prime $p$. Thus, there exists $x_0 \in \mathbb{Z}_p$ such that
$$x_0^2 - 2(8u+3)x_0 + N \equiv 0 \pmod{p}. \tag{31}$$
Since $N \not\equiv 0 \pmod{p}$, then $x_0 \neq 0$. Thus, from equation in (31) we have
$$u \equiv \frac{x_0^{-1}N + x_0 - 6}{16} \pmod{p} \tag{32}$$
for some $x_0 \in \mathbb{Z}_p \setminus \{0\}$.

Conversely, if there exists $u \in \left[\frac{-3+\sqrt{N}}{8}, \frac{N-15}{40}\right)$ where $(8u+3) - N$ is a perfect square, then
$$N = \left(8u + 3 + \sqrt{(8u+3)^2 - N}\right)\left(8u + 3 - \sqrt{(8u+3)^2 - N}\right). \tag{33}$$
If $N$ is prime, then $8u + 3 - \sqrt{(8u+1)^2 - N} = 1$, and hence, $u = m^2 + m$ which is a contradiction since $u \in \left[\frac{-3+\sqrt{N}}{8}, \frac{N-15}{40}\right)$. Therefore, $N$ is composite. □

The equation in (22) is called the Fermat's factorization of $16(2m+1)^2 + 1$. Note that statements 1-4 of proposition 2.4 states the property of all factors of $16(2m+1)^2 + 1$ if we factor $16(2m+1)^2 + 1$ using Fermat's factorization method. This also gives a sieve method to determine the proper factor of composite $16(2m+1)^2 + 1$. The next result characterized all proper factor of composite $16(2m+1)^2 + 1$ in terms of its structure.

**Proposition 2.5.** Every proper factors of composite $N = 16m^2 + 1$ is can be expressed as

$$8u + 3 \pm \sqrt{(8u+3)^2 - N} \tag{34}$$

where $u \in \mathbb{N}$ and
1. $(8u+3)^2 - N$ is a perfect square.
2. $u \in \left[\frac{-3+\sqrt{N}}{8}, \frac{N-15}{40}\right)$.
3. $4u + 1 \not\equiv 0 \pmod{p}$ for all prime $p \equiv 3 \pmod{4}$.
4. For all odd prime $p$ does not divide $N$, $u \equiv \frac{x_0^{-1}N + x_0 - 6}{16} \pmod{p}$, for some $x_0 \in \mathbb{Z}_p \setminus \{0\}$.

**Proof:** The results follow directly from Proposition 2.4. Furthermore, (34) follows from the fact that $8u + 3 \pm \sqrt{(8u+3)^2 - N}$ are roots of $f(x) = x^2 - 2(8u+3)x + N$. □

The next theorem is the main result of this study. The main result summarize the results in Proposition 2.3 and Proposition 2.5.

**Theorem 2.6.** Let $n \in \mathbb{N}$ and $N = 4n^2 + 1$.
1. If $n$ is even, then every proper factors of composite $N = 4n^2 + 1$ is can be expressed as
$$8u + 1 \pm \sqrt{(8u+1)^2 - N} \tag{35}$$
where $u \in \mathbb{N}$ and
1.1. $(8u + 1) - N$ is a perfect square.
1.2. $u \in \left[\frac{-1+\sqrt{N}}{8}, \frac{N-5}{40}\right)$.
1.3. $u \not\equiv 0 \pmod{p}$ for all prime $p \equiv 3 \pmod{4}$.
1.4. For all odd prime $p$ does not divide $N$, $u \equiv \frac{x_0^{-1}N + x_0 - 2}{16} \pmod{p}$, for some $x_0 \in \mathbb{Z}_p \setminus \{0\}$.

2. If $n$ is odd, then every proper factors of composite $N = 4n^2 + 1$ is can be expressed as
$$8u + 3 \pm \sqrt{(8u+3)^2 - N} \tag{36}$$
where $u \in \mathbb{N}$ and
2.1. $(8u+3)^2 - N$ is a perfect square.
2.2. $u \in \left[\frac{-3+\sqrt{N}}{8}, \frac{N-15}{40}\right)$.
2.3. $4u + 1 \not\equiv 0 \pmod{p}$ for all prime $p \equiv 3 \pmod{4}$.
2.4. For all odd prime $p$ does not divide $N$, $u \equiv \frac{x_0^{-1}N + x_0 - 6}{16} \pmod{p}$, for some $x_0 \in \mathbb{Z}_p \setminus \{0\}$.

**Proof:** Follows from Proposition 2.3 and Proposition 2.5. □

### III. FERMAT'S NUMBER

Fermat's number is a natural number of the form $F_n = 2^{2^n} + 1$ where $n$ is a nonnegative integer [2,10,11]. Note that determining the proper factors of composite Fermat's number is not easy by handful computation [2], for instance see [12, 13, 14, 15, 16, 17, 18, 19, 20, 21, 22, 23, 24, 25, 26]. In this section, we apply the same technique in section 2 to study the structure of proper factors of a Fermat's number.

Let $n$ be nonnegative integer and $p$ be prime factor of $F_n$. Lucas proved that if $n$ is nonnegative integer and $p$ is a prime factor of the Fermat's number $F_n = 2^{2^n} + 1$, then there exists a natural number $k$ where $p = 2^{n+2}k + 1$.[3] Consequently, every proper factor of composite Fermat's number is of the form $2^{n+2}k + 1$ where $k$ is a positive integer. Thus, the following lemma holds.

**Lemma 3.1.** Let $n \geq 4$. If the Fermat's number $F_n = 2^{2^n} + 1$ is composite, then there exists a natural number $s < \frac{\sqrt{F_n}-1}{2^{n+2}}$ where
$$2^{2^n - 2(n+2)} + s^2 \equiv 0 \ (mod(\ 2^{n+2}s + 1)). \tag{37}$$
Furthermore, $2^{n+2}s + 1$ is a proper factor of $F_n$.

**Proof:** Assume $F_n$ is composite. Then, by [3] there exists $r$ and $s$ such that:
$$2^{2^n} + 1 = (2^{n+2}s + 1)(2^{n+2}r + 1) \tag{38}$$
where $2^{n+2}s + 1 < F_n$ is a proper factor of $F_n$. Manipulate, then we have,
$$2^{2^n - (n+2)} = 2^{n+2}rs + r + s. \tag{39}$$
Since $n \geq 4$, then $2^n - 2(n+2) \geq 0$, and hence, there exists $\lambda \in \mathbb{N}$ such that $2^{n+2}\lambda = r + s$. Thus,
$$2^{2^n - (n+2)} = 2^{n+2}(2^{n+2}\lambda - s)s + 2^{n+2}\lambda \tag{40}$$
which gives
$$2^{2^n - 2(n+2)} + s^2 = \lambda(2^{n+2}s + 1) \tag{41}$$
where $s < \frac{\sqrt{F_n}-1}{2^{n+2}}$, as desired. □

The following proposition follows from lemma 3.1.

**Proposition 3.2.** Let $n \geq 5$. If the Fermat's number $F_n = 2^{2^n} + 1$ is composite, then there exists a natural number $\lambda$ where
1. $(2^{2n+3}\lambda + 1)^2 - F_n$ is a perfect square,
2. $\lambda \in \left[\frac{-1+\sqrt{F_n}}{2^{2n+3}}, 2^{2^n-(3n+5)}\right)$.
3. $\lambda \not\equiv 0 (mod\ p)$ where $p \equiv 3 (mod\ 4)$.
4. $\lambda \not\equiv 2 (mod 4)$.
5. $\lambda \equiv 1 (mod 3)$.

**Proof:**
1. Assume $F_n$ is composite. By lemma 3.1 there exists a natural number $s$ and integer $\lambda$ where
$$\frac{2^{2^n - 2(n+2)} + s^2}{2^{n+2}s + 1} = \lambda. \tag{42}$$
Since $s \in \mathbb{N}$, then $\lambda > 0$. Consider the quadratic polynomial
$$f(x) = x^2 - 2(2^{2n+3}\lambda + 1)x + F_n. \tag{43}$$
Since $f(2^{n+2}s + 1) = 0$ and the product of the roots of $f(x)$ is $F_n$, then every roots of $f(x)$ is a proper factor of $F_n$. Computing the zeroes $r_1$ and $r_2$ of $f(x)$, then we have
$$r_1, r_2 = 2^{2n+3}\lambda + 1 \pm \sqrt{(2^{2n+3}\lambda + 1)^2 - F_n}. \tag{44}$$
Thus, $\sqrt{(2^{2n+3}\lambda + 1)^2 - F_n} \in \mathbb{N} \cup \{0\}$, that is, $(2^{2n+3}\lambda + 1)^2 - F_n$ is a perfect square.

2. Note that we have
$$(2^{2n+3}\lambda + 1)^2 - F_n \geq 0. \tag{45}$$
Solving the inequality above with the assumption that $\lambda > 0$, then we have $\lambda \geq \frac{-1+\sqrt{F_n}}{2^{2n+3}}$. In addition, by Rolle's theorem, there exists $\theta \in [r_1, r_2] \subset [r_1, 2^{2^n-(n+2)} + 1]$ where $f'(\theta) = 0$. Thus,
$$2\theta - 2(2^{2n+3}\lambda + 1) = 0. \tag{46}$$
Thus, $\lambda = \frac{\theta - 1}{2^{2n+3}}$. Since $\theta \in [r_1, r_2] \subset [r_1, 2^{2^n-(n+2)} + 1]$, then $\lambda \leq 2^{2^n-(3n+5)}$. Therefore, we have $\lambda \in \left[\frac{-1+\sqrt{F_n}}{2^{2n+3}}, 2^{2^n-(3n+5)}\right)$.
3. Follows from statement 3 of proposition 2.2.
4. Follows from statement 5 of proposition 2.2.
5. Follows from statement 6 of proposition 2.2. □

The following theorem are direct consequence of proposition 3.2.

**Theorem 3.3.** Let $n \geq 5$. Then, every proper factors of $F_n$ is can be expressed as
$$2^{2n+3}\lambda + 1 \pm \sqrt{(2^{2n+3}\lambda + 1)^2 - F_n} \tag{47}$$

where $\lambda \in \mathbb{N}$ and
1. $(2^{2n+3}\lambda + 1)^2 - F_n$ is a perfect square
2. $\lambda \in \left[\frac{-1+\sqrt{F_n}}{2^{2n+3}}, 2^{2^n-(3n+5)}\right)$
3. $\lambda \not\equiv 0 \pmod{p}$ where $p \equiv 3 \pmod{4}$
4. $\lambda \not\equiv 2 \pmod{4}$
5. $\lambda \equiv 1 \pmod{3}$

**Proof:** Theorem 3.3 follows from proposition 3.2. Furthermore, equation (47) follows from (44). □

## VI. CONCLUSIONS

In this study, we give a characterization of all proper factors of $4n^2 + 1$ in terms of its form by applying Fermat's method. In addition, we derived the property of all factors of $4n^2 + 1$ that depends on the parity of $n$ (see Theorem 2.6). In addition, the results in Proposition 2.2 and Proposition 2.5 give a new sieve method to determine the factors of $4n^2 + 1$. Fermat's number is also considered in this study by deriving a new property of composite Fermat's number that is similar in proposition 2.2 and theorem 2.6 (see Proposition 3.2 and Theorem 3.2).

## ACKNOWLEDGMENT
The authors are grateful to the Department of Mathematics and Statistics of Polytechnic University of the Philippines, Manila for their unending support to finish this paper.